\documentclass[11pt]{amsart}
\usepackage{xcolor}

\include{package}
\definecolor{grn}{rgb}{0,0.6,0}
\definecolor{mrn}{rgb}{0.3,0,0}
\definecolor{blue}{rgb}{0,0,0.7}
\definecolor{Mygray}{rgb}{0.75,0.75,0.75}
\definecolor{auburn}{rgb}{0.43, 0.21, 0.1}
\definecolor{britishracinggreen}{rgb}{0.0, 0.26, 0.15}
\definecolor{taupe}{rgb}{0.28, 0.24, 0.2}
\usepackage{amsmath,amssymb}
\newtheorem{theorem}{Theorem}

\newtheorem{proposition}{Proposition}

\newtheorem{lemma}{Lemma}
\theoremstyle{definition}
\newtheorem{definition}{Definition}
\newtheorem{question}{Question}
\newtheorem{rmk}{Remark}

\usepackage{latexsym,amssymb}
\begin{document}
\baselineskip=14.5pt
\title[Non-P\'{o}lya bi-quadratic fields]{Non-P\'{o}lya bi-quadratic fields with an Euclidean ideal class} 

\author{Jaitra Chattopadhyay and Anupam Saikia}
\address[Jaitra Chattopadhyay]{Department of Mathematics, Indian Institute of Technology, Guwahati, Guwahati - 781039, Assam, India}
\email[Jaitra Chattopadhyay]{jaitra@iitg.ac.in}
\address[Anupam Saikia]{Department of Mathematics, Indian Institute of Technology, Guwahati, Guwahati - 781039, Assam, India}
\email[Anupam Saikia]{a.saikia@iitg.ac.in}

\begin{abstract}
For an integral domain $R$, the {\it ring of integer-valued polynomials} over $R$ consists of all polynomials $f(X) \in R[X]$ such that $f(R) \subseteq R$. An interesting case to study is when $R$ is a Dedekind domain, in particular when $R$ is the ring of integers of an algebraic number field. An algebraic number field $K$ with ring of integers $\mathcal{O}_{K}$ is said to be a P\'{o}lya field if the $\mathcal{O}_{K}$-module of integer-valued polynomials on $K$ admits a regular basis. Associated to $K$ is a subgroup $Po(K)$ of the ideal class group $Cl_{K}$, known as the {\it P\'{o}lya group of $K$}, that measures the failure of $K$ from being a P\'{o}lya field. In this paper, we prove the existence of three pairwise distinct totally real bi-quadratic fields, each having P\'{o}lya group isomorphic to $\mathbb{Z}/2\mathbb{Z}$. This extends the previously known families of number fields considered by Heidaryan and Rajaei in \cite{rajaei-jnt} and \cite{rajaei}. Our results also establish that under mild assumptions, the possibly infinite families of bi-quadratic fields having a non-principal Euclidean ideal class, considered in \cite{self-jnt}, fail to be P\'{o}lya fields.
\end{abstract}

\renewcommand{\thefootnote}{}



\renewcommand{\thefootnote}{\arabic{footnote}}
\setcounter{footnote}{0}

\maketitle

\section{Introduction}

Let $R$ be an integral domain with field of fractions $K$. The polynomial ring $K[X]$ is a central object of study in algebra as well as in algebraic number theory. In particular, the set of polynomials over $K$ that carry $R$ to itself, occupies a special place in studying various properties of $K[X]$. We define $$Int(R) = \{f(X) \in K[X] : f(R) \subseteq R\}.$$ It is a trivial observation that $R[X] \subseteq Int(R) \subseteq K[X]$. Moreover, $Int(R)$ is an additive subgroup of $K[X]$ and there is a natural $R$-module structure on $Int(R)$, defined by $(\alpha, f(X)) \mapsto \alpha \cdot f(X)$, where $\alpha \in R$ and $f(X) \in Int(R)$. 



\smallskip

Several interesting properties emerge for particular choices of the ring $R$. Among those, one of the most important cases is when $K$ is an algebraic number field and $R$ is its ring of integers. In that case, the object $Int(\mathcal{O}_{K})$ has deep connections with the ideal class group $Cl_{K}$ of $K$. As the class number $h_{K}$ of $K$ measures the failure of unique factorization of elements of $\mathcal{O}_{K}$ into the product of irreducible elements, likewise we associate a subgroup of $Cl_{K}$ that measures the extent of failure of the $\mathcal{O}_{K}$-module $Int(\mathcal{O}_{K})$ from being a nicely behaved $\mathcal{O}_{K}$-module. Before we introduce that particular subgroup, called ``P\'{o}lya group", of $Cl_{K}$, we begin with the definition of a regular basis which motivates the study of P\'{o}lya groups and $Int(\mathcal{O}_{K})$.

\begin{definition} (cf. \cite{cahen-chabert-book})
A basis $\{f_{j}\}_{j \geq 0}$ of the $\mathcal{O}_{K}$-module $Int(\mathcal{O}_{K})$ is said to be {\it regular} if degree$(f_{j}) = j$ for all $j \geq 0$.
\end{definition}

Using the notion of a regular basis, {\it P\'{o}lya fields} are defined as follows.

\begin{definition} (cf. \cite{polya})\label{defn-polya-field}
A number field $K$ is said to be a P\'{o}lya field if the $\mathcal{O}_{K}$-module $Int(\mathcal{O}_{K})$ admits a regular basis.
\end{definition}

It is easy to note that the polynomials $\{f_{n}\}_{n \geq 0}$, where $f_{n} = {{X}\choose{n}} = \frac{X(X - 1)\ldots (X - n +1)}{n!}$ for $n \geq 1$ and $f_{0} = 1$, constitute a regular basis for $Int(\mathbb{Z})$. Consequently, $\mathbb{Q}$ is a P\'{o}lya field. But an arbitrary number field $K$ may not have a regular basis. Therefore, it is useful to look for conditions on the number field $K$ which imply that $K$ is indeed a P\'{o}lya field. For that, we first introduce the {\it characteristic ideals} $\mathfrak{J}_{n}(K)$ of $K$ as follows.

\begin{definition} (cf. \cite{cahen-chabert-book})
For an integer $n \geq 1$, let $$\mathfrak{J}_{n}(K) = \{\mbox{leading coefficients of elements of } Int(\mathcal{O}_{K}) \mbox{ of degree } n\} \cup \{0\}.$$
\end{definition}

It can be checked that for each integer $n \geq 1$, the set $\mathfrak{J}_{n}(K)$ is a fractional ideal of $K$. The following proposition provides us a connection between $\mathfrak{J}_{n}(K)$ and $K$ being a P\'{o}lya field.

\begin{proposition} (cf. \cite{cahen-chabert-book})\label{char-th}
A number field $K$ is a P\'{o}lya field if and only if $\mathfrak{J}_{n}(K)$ is principal for all integer $n \geq 1$.
\end{proposition}

By Proposition \ref{char-th}, it is immediate that if $\mathcal{O}_{K}$ is a principal ideal domain, then $\mathfrak{J}_{n}(K)$ is also principal for all integer $n \geq 1$. In other words, if $h_{K} = 1$, then $K$ is a P\'{o}lya field. This provides us a class of number fields that are P\'{o}lya. But if $h_{K} \neq 1$, then we cannot directly appeal to Proposition \ref{char-th} to decide whether $K$ is a P\'{o}lya field or not. This motivates us to study the class group $Cl_{K}$ and its subgroups that might shed some light on the algebraic structure of $Int(\mathcal{O}_{K})$. For that, we now define a particular subgroup of $Cl_{K}$ that captures the extent of failure of $K$ from being a P\'{o}lya field.

\begin{definition} (cf. \cite{cahen-chabert-book})\label{polya-group}
Let $K$ be a number field with ring of integers $\mathcal{O}_{K}$ and ideal class group $Cl_{K}$. For each integer $n \geq 1$, let $[\mathfrak{J}_{n}(K)]$ be the ideal class in $Cl_{K}$ corresponding to the fractional ideal $\mathfrak{J}_{n}(K)$. The P\'{o}lya group $Po(K)$ of $K$ is defined to be the subgroup of $Cl_{K}$ generated by the elements $[\mathfrak{J}_{n}(K)]$ in $Cl_{K}$.
\end{definition}

We introduce another quantity which is also related to $K$ being a P\'{o}lya field. For an integer $q \geq 1$, we define $\Pi_{q}(K) = \displaystyle\prod_{\substack{\mathfrak{p} \mbox{ prime } \\ N(\mathfrak{p}) = q}} \mathfrak{p}$. If $q$ is not the norm of any prime ideal in $\mathcal{O}_{K}$, then we define $\Pi_{q}(K) = \mathcal{O}_{K}$. The following Proposition gives us a necessary and sufficient condition for $K$ to be a P\'{o}lya field in terms of the invariants introduced so far.

\begin{proposition} (cf. \cite{cahen-chabert-book}, \cite{ostrowski})\label{EQIVALENCE}
A number field $K$ is a P\'{o}lya field if and only if one of these three equivalent conditions holds.

\begin{enumerate}

\item $\mathfrak{J}_{n}(K)$ is principal for all integer $n \geq 1$.

\item $\Pi_{q}(K)$ is principal for each integer $q \geq 2$.

\item $Po(K) = \{1\}$.
\end{enumerate}
\end{proposition}

\begin{rmk}
It is evident from Proposition \ref{EQIVALENCE} that the question of studying the integer-valued polynomials on $K$ can be translated into a problem of the study of the product of prime ideals in $\mathcal{O}_{K}$ with a given norm. 
\end{rmk}

It is quite natural to look for explicit families of number fields that are P\'{o}lya. In fact, the classification of P\'{o}lya fields of lower degree number fields is of considerable interest. Towards that direction, a complete characterization of quadratic P\'{o}lya fields is known due to the work of Zantema \cite{zantema}.
\begin{proposition} \cite{zantema}\label{quadratic-poly-classification}
Let $d \neq 1$ be a square-free integer and let $K = \mathbb{Q}(\sqrt{d})$. Then $K$ is a P\'{o}lya field if and only if one of the following conditions holds.
\begin{enumerate}
\item $d = -1, -2, 2$,

\item $d = -p$, where $p \equiv 3 \pmod {4}$ is a prime number,

\item $d = p$, where $p > 0$ is an odd prime number,

\item $d = 2p$, where $p$ is a prime number and either $p \equiv 3 \pmod {4}$ or $p \equiv 1 \pmod {4}$ and the fundamental unit of $K$ has norm $1$,

\item $d = pq$, where $p$ and $q$ are prime numbers with either $p \equiv q \equiv 3 \pmod {4}$ or $p \equiv q \equiv 1 \pmod {4}$ and the fundamental unit of $K$ has norm $1$.
\end{enumerate}
\end{proposition}

In \cite{leriche}, Leriche classified cubic, quartic and sextic P\'{o}lya fields under certain assumptions. Moreover, in \cite{leriche-jtnb}, Leriche completely classified bi-quadratic P\'{o}lya fields that are compositum of two quadratic P\'{o}lya fields. More precisely, she proved the following.
\begin{proposition} \cite{leriche}\label{bi-quadratic-polya-classification}
Let $m \neq 1$ and $n \neq 1$ be two square-free integers such that both $\mathbb{Q}(\sqrt{m})$ and $\mathbb{Q}(\sqrt{n})$ are P\'{o}lya fields. The the bi-quadratic field $\mathbb{Q}(\sqrt{m}, \sqrt{n})$ is a P\'{o}lya field except for the following cases.
\begin{enumerate}
\item $\mathbb{Q}(\sqrt{-2}, \sqrt{p})$ is not a P\'{o}lya field, where $p$ is a prime with $p \equiv 3 \pmod {4}$.

\item $\mathbb{Q}(\sqrt{-1}, \sqrt{2q})$ is not a P\'{o}lya field, where $q$ is an odd prime number.

\end{enumerate}

Further, if for prime numbers $p$ and $q$, the field $\mathbb{Q}(\sqrt{p}, \sqrt{2q})$ is a P\'{o}lya field, then 
\begin{enumerate}
\item either $p \equiv -1 \pmod {8}$ and $q \equiv \pm {1} \pmod {8}$

\item or $p \equiv 3 \pmod {8}$ and $q \equiv 1, 3 \pmod {8}$.
\end{enumerate}

\end{proposition}

In view of Proposition \ref{bi-quadratic-polya-classification}, it is interesting to ask for the P\'{o}lya-ness of a bi-quadratic field that is not a compositum of two quadratic P\'{o}lya fields. Heidaryan and Rajaei pursued this question in \cite{rajaei-jnt} and \cite{rajaei} and obtained some such families of number fields. In the next section, we state our main results of this paper, which are extensions of the families of bi-quadratic fields considered in \cite{rajaei-jnt} and \cite{rajaei}.

\section{Statements of main theorems}

\begin{theorem}\label{main-1} 
Let $p,q$ and $r$ be distinct prime numbers with $p \equiv 3 \pmod {4}$ and $q \equiv r \equiv 1 \pmod {8}$. Assume that $\left(\frac{q}{r}\right) = -1$, where $\left(\frac{\cdot}{r}\right)$ stands for the Legendre symbol. Then for the bi-quadratic field $K^{\prime} = \mathbb{Q}(\sqrt{p},\sqrt{qr})$, we have $Po(K^{\prime}) \simeq \mathbb{Z}/2\mathbb{Z}$. In particular, $K^{\prime}$ is not a P\'{o}lya field.
\end{theorem}


\begin{theorem}\label{NEW-TH}
Let $p,q$ and $r$ be distinct prime numbers with $p \equiv q \equiv 3 \pmod {4}$ and $r \equiv 1 \pmod {8}$. Assume that $\left(\frac{p}{r}\right) = 1$ and $\left(\frac{q}{r}\right) = -1$. Then for the the bi-quadratic field $F = \mathbb{Q}(\sqrt{p},\sqrt{qr})$, we have $Po(F) \simeq \mathbb{Z}/2\mathbb{Z}$. In particular, $F$ is not a P\'{o}lya field.
\end{theorem}


\begin{theorem}\label{main-2}
Let $p$ and $q$ be distinct prime numbers with $p \equiv q \equiv 1 \pmod {4}$. Assume that $\left(\frac{p}{q}\right) = -1$. Then for the bi-quadratic field $K^{\prime \prime} = \mathbb{Q}(\sqrt{2},\sqrt{pq})$, we have $Po(K^{\prime \prime}) \simeq \mathbb{Z}/2\mathbb{Z}$. In particular, $K^{\prime \prime}$ is not a P\'{o}lya field.
\end{theorem}


\begin{rmk}\label{LATEST REMARK}
In \cite{rajaei}, it was proved that the totally real bi-quadratic field $\mathbb{Q}(\sqrt{p},\sqrt{qr})$ is a P\'{o}lya field whenever $p,q$ and $r$ are distinct prime numbers with $p \equiv q \equiv 3 \pmod {4}$ and $r \equiv 5 \pmod {8}$. Theorem \ref{NEW-TH} shows that for $r \equiv 1 \pmod {8}$, the field $\mathbb{Q}(\sqrt{p},\sqrt{qr})$ is not P\'{o}lya under the mild assumptions on the Legendre symbols.
\end{rmk}


\begin{rmk}
A result of Dirichlet (cf. \cite{trotter}) implies that for two distinct prime numbers $r$ and $s$, if the Legendre symbol $\left(\frac{r}{s}\right) = -1$, then the fundamental unit of $\mathbb{Q}(\sqrt{rs})$ has norm $-1$. Thus the hypotheses of Theorem \ref{main-1}, Theorem \ref{NEW-TH} and Theorem \ref{main-2} imply that the fundamental units of the subfields $\mathbb{Q}(\sqrt{qr})$ and $\mathbb{Q}(\sqrt{pq})$ of $K^{\prime}$, $F$ and $K^{\prime\prime}$ respectively, have norm $-1$. Consequently, by Proposition \ref{quadratic-poly-classification}, the quadratic fields $\mathbb{Q}(\sqrt{qr})$ and $\mathbb{Q}(\sqrt{pq})$ are not P\'{o}lya. Thus $K^{\prime}$, $F$ and $K^{\prime\prime}$ each contains precisely one quadratic subfield which is P\'{o}lya.
\end{rmk}


\begin{rmk}
It is interesting to note that the number fields considered in Theorem \ref{main-1}, Theorem \ref{NEW-TH} and Theorem \ref{main-2} have even class numbers because the P\'{o}lya group of a number field is a subgroup of the ideal class group and in each of the three above theorems, the P\'{o}lya group is asserted to be $\mathbb{Z}/2\mathbb{Z}$.
\end{rmk}

\begin{rmk}
The conditions on the Legendre symbols on the primes $p,q$ and $r$ are very crucial. In fact, a result of Pollack \cite{paul} implies that for a given prime number $r \geq 13$, there exist prime numbers $p \equiv 3 \pmod {4}$ and $q \equiv 1 \pmod {4}$ such that $\left(\frac{p}{r}\right) = \left(\frac{q}{r}\right) = -1$ and $2 \leq p,q \leq r - 1$. Even though there are plenty of prime numbers satisfying the hypotheses of Theorem \ref{main-1}, Theorem \ref{NEW-TH} and Theorem \ref{main-2}, the result of Pollack guarantees that the hypotheses hold for small primes as well.
\end{rmk}

\section{Connection with Euclidean ideal class}\label{subsection}
We have earlier observed that if the class group $Cl_{K}$ is trivial, then $K$ is necessarily a P\'{o}lya field. It is natural to investigate the P\'{o}lya-ness of $K$ under various other assumptions on $Cl_{K}$. One such instance is when $Cl_{K}$ is a cyclic group. Interestingly, the group $Cl_{K}$ being cyclic is closely related to the existence of a non-principal {\it Euclidean ideal class} in $K$. First, following \cite{lenstra}, we define an Euclidean ideal as follows.
\begin{definition} \cite{lenstra}
Let $K$ be a number field with ring of integers $\mathcal{O}_{K}$. Let $\mathbb{I}$ be the set of all fractional ideals containing $\mathcal{O}_{K}$ and let $W$ be a well-ordered set. A fractional ideal $C$ of $\mathcal{O}_{K}$ is said to be an Euclidean ideal if there exists a function $\psi : \mathbb{I} \rightarrow W$ such that for any $J \in \mathbb{I}$ and any $x \in JC \setminus C$, there exists some $\alpha \in C$ satisfying 
\begin{equation}
\psi((x-\alpha)^{-1}JC) < \psi (J).
\end{equation}
\end{definition}

Lenstra \cite{lenstra} proved that if $C$ is an Euclidean ideal, then every fractional ideal $C^{\prime}$ belonging to the ideal class $[C]$ is also an Euclidean ideal. Therefore, an ideal class $[C]$ is said to be an {\it Euclidean ideal class}, abbreviated as EIC, if $C$ is an Euclidean ideal. One of the most important consequences of the existence of an EIC in $K$ is that $Cl_{K}$ is cyclic. Lenstra \cite{lenstra} proved under the assumption of the Extended Riemann Hypothesis (ERH) that the converse holds true for all number fields $K$ with ${\rm{rk}}(\mathcal{O}_{K}^{*}) \geq 1$. Later, many explicit families of number fields with cyclic class groups were unconditionally proved to possess an EIC (cf. \cite{self-jnt}, \cite{graves-ijnt}, \cite{cathy}).

\smallskip

Thus it is equally interesting to inquire the nature of the P\'{o}lya group of a number field having an EIC. In \cite{self-jnt}, it has been proved that for prime numbers $p \equiv 3 \pmod {4} \mbox{ and } q \equiv r \equiv 1 \pmod {4}$, the totally real bi-quadratic fields $K_{1} = \mathbb{Q}(\sqrt{p},\sqrt{qr})$ and $K_{2} = \mathbb{Q}(\sqrt{2},\sqrt{qr})$ have an EIC provided $h_{K_{1}} = h_{K_{2}} = 2$. Therefore, by Theorem \ref{main-1} and Theorem \ref{main-2}, we can conclude that there exist number fields with an EIC even though they fail to be P\'{o}lya.


\section{Preliminaries}

In this section, we record the necessary results needed to prove Theorem \ref{main-1} and Theorem \ref{main-2}. The main tool in proving our theorems is to make use of a result due to Zantema \cite{zantema}. We briefly furnish the result as used in \cite{zantema}.

\smallskip

Let $K$ be a finite Galois extension of $\mathbb{Q}$ with Galois group $Gal(K/\mathbb{Q}) = G$. Then there is a canonical action of $G$ on the group of units $\mathcal{O}_{K}^{*}$ defined by $(\sigma,\alpha) \mapsto \sigma (\alpha)$, where $\sigma \in G$ and $\alpha \in \mathcal{O}_{K}^{*}$. Thus, $\mathcal{O}_{K}^{*}$ is endowed with a $G$-module structure. Zantema's result connects the first cohomology group $H^{1}(G,\mathcal{O}_{K}^{*})$ with the ramified primes in $K/\mathbb{Q}$ as follows.

\begin{proposition} \cite{zantema}\label{zantema's-main-theorem}
Let $K$ be a finite Galois extension of $\mathbb{Q}$ with Galois group $Gal(K/\mathbb{Q}) = G$. For a rational prime $\ell$, let $e_{\ell}$ be the ramification index of $\ell$ in $K/\mathbb{Q}$. Then there exist a canonical embedding $\Psi : H^{1}(G,\mathcal{O}_{K}^{*}) \to \displaystyle\bigoplus_{\ell \mbox{ prime }}\mathbb{Z}/e_{\ell}\mathbb{Z}$ and an exact sequence of abelian groups $$1 \rightarrow H^{1}(G,\mathcal{O}_{K}^{*}) \xrightarrow{\Psi} \displaystyle\bigoplus_{\ell \mbox{ prime}} \mathbb{Z}/e_{\ell} \mathbb{Z} \rightarrow Po(K) \rightarrow 1.$$ 
\end{proposition}

As an immediate corollary to Proposition \ref{zantema's-main-theorem}, we see that $K$ is a P\'{o}lya field if and only if $|H^{1}(G,\mathcal{O}_{K}^{*})| = \displaystyle\prod_{\ell \mbox{ ramified}}e_{\ell}$.

\smallskip

Next, we record two lemmas that will be useful to understand the group $H^{1}(G,\mathcal{O}_{K}^{*})$ in a more convenient way when $K$ a bi-quadratic field.

\begin{lemma} \cite{bennet}\label{lem-1}
Let $K$ be a bi-quadratic field with quadratic subfields $K_{1}, K_{2}$ and $K_{3}$. Let $H$ be the subgroup of $H^{1}(G,\mathcal{O}_{K}^{*})$ consisting of elements of order dividing $2$. Then the index of $H$ in $H^{1}(G,\mathcal{O}_{K}^{*})$ is $\leq 2$. The index is $2$ if and only if the rational prime $2$ is totally ramified in $K/\mathbb{Q}$ and there exist $\alpha_{i} \in \mathcal{O}_{K_{i}}$ for $i = 1,2,3$ such that $$N(\alpha_{1}) = N(\alpha_{2}) = N(\alpha_{3}) = \pm {2}.$$
\end{lemma}

\smallskip

\begin{lemma} \cite{zantema}\label{lem-2}
Let $K$ be a bi-quadratic field with quadratic subfields $K_{1}, K_{2}$ and $K_{3}$. Let $\Delta_{i}$ be the square-free part of the discriminant $d_{K_{i}}$ of $K_{i}$, for $i = 1,2,3$. Let $u_{i} = z_{i} + t_{i}\sqrt{\Delta_{i}}$ be a fundamental unit of $\mathcal{O}_{K_{i}}$ with $z_{i} > 0$. Define 
\begin{equation*}
a_{i}=
\begin{cases}
    N(u_{i} + 1)   & ~ \text{ if } N(u_{i}) = 1,\\
    1 &  ~ \text{ if } N(u_{i}) = -1.
\end{cases}
\end{equation*}
Let $H$ be the subgroup of $H^{1}(G,\mathcal{O}_{K}^{*})$ consisting of elements of order dividing $2$. Then $H$ is generated by the images of $\Delta_{1}, \Delta_{2}, \Delta_{3}, a_{1}, a_{2}$ and $a_{3}$ in $\mathbb{Q}^{*}/(\mathbb{Q}^{*})^{2}$.
\end{lemma}

\smallskip

Lemma \ref{lem-1} and Lemma \ref{lem-2} have been used in \cite{rajaei-jnt} and \cite{rajaei} to prove their main theorems. We shall adapt the similar line of approach as taken in \cite{rajaei} to prove Theorem \ref{main-1} and Theorem \ref{main-2}. Before we go to the proofs of our main theorems, let us state the following lemma from \cite{rajaei}, without proof.

\begin{lemma} \cite{rajaei}\label{lem-0}
Let $p$ be a prime number with $p \equiv 3 \pmod {4}$. Let $u = z + t\sqrt{p}$ be a fundamental unit of $\mathbb{Q}(\sqrt{p})$. Then $1+z$ and $t$ are odd integers. Moreover, if $p \equiv 7 \pmod {8}$, then $1+z$ is a square and if $p \equiv 3 \pmod {8}$, then $1+z = pm^{2}$ for some integer $m$.
\end{lemma}

\section{Proof of theorem \ref{main-1}}
Let $K_{1} = \mathbb{Q}(\sqrt{p}), K_{2} = \mathbb{Q}(\sqrt{qr})$ and $K_{3} = \mathbb{Q}(\sqrt{pqr})$. Since $p \equiv 3 \pmod {4}$ and $q \equiv r \equiv 1 \pmod {8}$, we obtain that precisely the rational primes $2,p,q$ and $r$ are ramified in $K^{\prime}/\mathbb{Q}$. If $e_{\ell}$ stands for the ramification index of a rational prime $\ell$ in $K^{\prime}/\mathbb{Q}$, then we have $$e_{2} = e_{p} = e_{q} = e_{r} = 2.$$ Therefore, by Proposition \ref{zantema's-main-theorem}, we have that $H^{1}(G,\mathcal{O}_{K^{\prime}}^{*})$ embeds into $\mathbb{Z}/2\mathbb{Z} \oplus \mathbb{Z}/2\mathbb{Z} \oplus \mathbb{Z}/2\mathbb{Z} \oplus \mathbb{Z}/2\mathbb{Z}$. Moreover, if we prove that $|H^{1}(G,\mathcal{O}_{K^{\prime}})| = 8$, then it follows from Proposition \ref{zantema's-main-theorem} that $$Po(K^{\prime}) \simeq \displaystyle\bigoplus_{\ell \mbox{ prime}} \mathbb{Z}/e_{\ell} \mathbb{Z}\Bigg/H^{1}(G,\mathcal{O}_{K^{\prime}}^{*}) \simeq \mathbb{Z}/2\mathbb{Z}.$$ This is why, it is enough to prove that $|H^{1}(G,\mathcal{O}^{*}_{K^{\prime}})| = 8$.

\smallskip

We note that the rational prime $2$ is not totally ramified in $K^{\prime}/\mathbb{Q}$ because $2$ is unramified in $K_{2}/\mathbb{Q}$. Consequently, by Lemma \ref{lem-1}, we conclude that $H = H^{1}(G,\mathcal{O}^{*}_{K^{\prime}})$. We prove that $|H| = 8$.

\smallskip

In the notations of Lemma \ref{lem-2}, we have $\Delta_{1} = p$, $\Delta_{2} = qr$ and $\Delta_{3} = pqr$. We also use the notation $[\alpha]$ to denote the image of the element $\alpha\in \mathbb{Q}$ in $\mathbb{Q}^{*}/(\mathbb{Q}^{*})^{2}$. Hence $[\Delta_{1}], [\Delta_{2}], [\Delta_{3}] \in \langle [p], [qr] \rangle$ in $\mathbb{Q}^{*}/(\mathbb{Q}^{*})^{2}$. For $i \in \{1,2,3\}$, let $u_{i} = z_{i} + t_{i}\sqrt{\Delta_{i}}$ be a fundamental unit of $K_{i}$. We have $N(u_{1}) = z_{1}^{2} - t_{1}^{2}p$. If $N(u_{1}) = -1$, then reducing the equation $\pmod {p}$, we get $z_{1}^{2} \equiv -1\pmod {p}$, which is impossible because $p \equiv 3 \pmod {4}$. Therefore, $N(u_{1}) = 1$ and hence $a_{1} = 2(z_{1} + 1)$. By Lemma \ref{lem-0}, we have $[a_{1}] = [2] \mbox{ or } [2p]$. 

\smallskip

Since $\left(\frac{q}{r}\right) = -1$, by a result of Dirichlet (cf. \cite{trotter}), we get that the fundamental unit of $K_{2} = \mathbb{Q}(\sqrt{qr})$ has norm $-1$. Therefore, by Lemma \ref{lem-2}, we have $a_{2} = 1$.

\smallskip

Let $u = z + t\sqrt{pqr}$ be a fundamental unit of $K_{3} = \mathbb{Q}(\sqrt{pqr})$. Again, the condition $p \equiv 3 \pmod {4}$ implies that $N(u) = 1$. That is, $z^{2} - t^{2}pqr = 1$ and hence $(z-1)(z+1) = t^{2}pqr$.

\smallskip

\noindent
{\bf Case 1.} $z$ is an even integer. Therefore, $t$ is an odd integer. Moreover, $z-1$ and $z+1$ are two consecutive odd integers and therefore $\gcd(z-1,z+1) = 1$. This, together with the equation $(z-1)(z+1) = t^{2}pqr$, implies that $z - 1 = n^{2}\eta$ and $z + 1 = m^{2}\epsilon$ for coprime odd integers $m, n$ with $mn = t$ and coprime integers $\epsilon, \eta$ with $\epsilon\eta = pqr$. Consequently, we have $$[a_{3}] = [2(z + 1)] = [2m^{2}\epsilon] = [2\epsilon] \in \langle [2],[p],[qr] \rangle \subseteq \mathbb{Q}^{*}/(\mathbb{Q}^{*})^{2} \mbox{ if and only if } \epsilon = 1,p,qr,pqr.$$ We now prove that $\epsilon = q,r,pq \mbox{ and } pr$ cannot occur.

\smallskip

First, we note that eliminating $z$ from the equations $z + 1 = m^{2}\epsilon$ and $z - 1 = n^{2}\eta$, we obtain $2 = m^{2}\epsilon - n^{2}\eta$. If $\epsilon = q$, then $\eta = pr$ and therefore $2 = m^{2}q - n^{2}pr$. Hence we have $1 = \left(\frac{2}{r}\right) = \left(\frac{q}{r}\right) = -1$, which is a contradiction. Again, if $\epsilon = r$, then $\eta = pq$ and $2 = m^{2}r - n^{2}pq$. Consequently, we have $1 = \left(\frac{2}{q}\right) = \left(\frac{r}{q}\right) = -1$, which is a contradiction. Further, if $\epsilon = pq$, then $\eta = r$ and therefore $2 = m^{2}pq - n^{2}r$. Hence we have $1 = \left(\frac{2}{q}\right) = \left(\frac{-r}{q}\right) = \left(\frac{r}{q}\right) = -1$, which is a contradiction. Lastly, if $\epsilon = pr$, then $\eta = q$ and therefore $2 = m^{2}pr - n^{2}q$. Hence we have $1 = \left(\frac{2}{r}\right) = \left(\frac{-q}{r}\right) = \left(\frac{q}{r}\right) = -1$, which yields a contradiction.

\medskip

\noindent
{\bf Case 2.} $z$ is an odd integer. In this case, both $z - 1$ and $z + 1$ are even integers and therefore, $\gcd(z - 1, z+ 1) = 2$. From the equation $(z - 1)(z + 1) = t^{2}pqr$, we see that $t$ is an even integer. Therefore, $t = 2k$ for some integer $k \geq 1$. Thus $\frac{z - 1}{2}\cdot \frac{z + 1}{2} = k^{2}pqr$ and $\gcd\left(\frac{z - 1}{2},\frac{z + 1}{2}\right) = 1$. Hence we have $\frac{z - 1}{2} = n^{2}\eta$ and $\frac{z + 1}{2} = m^{2}\epsilon$ for coprime integers $m,n$ with $mn = k$ and coprime integers $\epsilon,\eta$ with $\epsilon\eta = pqr$. Therefore, $$[a_{3}] = [2(z + 1)] = [4m^{2}\epsilon] = [\epsilon] \in \langle [2],[p],[qr] \rangle \mbox{ if and only if } \epsilon = 1,p,qr,pqr.$$ We now prove that $\epsilon = q,r,pq \mbox{ and } pr$ cannot occur. For that, we first note that $\frac{z - 1}{2} = n^{2}\eta$ and $\frac{z + 1}{2} = m^{2}\epsilon$ yield $1 = m^{2}\epsilon - n^{2}\eta$ .

\smallskip

If $\epsilon = q$, then $\eta = pr$ and therefore $1 = m^{2}q - n^{2}pr$. Hence $1 = \left(\frac{1}{r}\right) = \left(\frac{m^{2}q}{r}\right) = \left(\frac{q}{r}\right) = -1$, which is a contradiction. Again, if $\epsilon = r$, then $\eta = pq$ and $1 = m^{2}r - n^{2}pq$. Hence $1 = \left(\frac{1}{q}\right) = \left(\frac{m^{2}r}{q}\right) = \left(\frac{r}{q}\right) = -1$, which is a contradiction. Similarly, if $\epsilon = pq$, then $\eta = r$ and we have $1 = m^{2}pq - n^{2}r$. Hence $1 = \left(\frac{1}{q}\right) =  \left(\frac{-r}{q}\right) = \left(\frac{r}{q}\right) = -1$, which yields a contradiction. Lastly, if $\epsilon = pr$, then $\eta = q$ and we have $1 = m^{2}pr - n^{2}q$. Hence $1 = \left(\frac{1}{r}\right) = \left(\frac{-q}{r}\right) = \left(\frac{q}{r}\right) = -1$, which is a contradiction.

\smallskip

Thus we obtain that $$\langle [\Delta_{1}], [\Delta_{2}], [\Delta_{3}], [a_{1}], [a_{2}], [a_{3}]\rangle = \langle [2], [p], [qr] \rangle \subseteq \mathbb{Q}^{*}/(\mathbb{Q}^{*})^{2}.$$ 

Consequently, we have $|H^{1}(G,\mathcal{O}^{*}_{K^{\prime}})| = |H| = 8$ and hence $Po(K^{\prime}) \simeq \mathbb{Z}/2\mathbb{Z}$. This completes the proof of Theorem \ref{main-1}. $\hfill\Box$

\medskip

\section{Proof of Theorem \ref{NEW-TH}}

Let $F_{1} = \mathbb{Q}(\sqrt{p})$, $F_{2} = \mathbb{Q}(\sqrt{qr})$ and $F_{3} = \mathbb{Q}(\sqrt{pqr})$. The rational primes that are ramified in $F/\mathbb{Q}$ are precisely $2,p,q$ and $r$. Also, the ramification indices of $p,q$ and $r$ are $2$ each since $p$ is unramified in $F_{2}/\mathbb{Q}$ and $q$ and $r$ are unramified in $F_{1}/\mathbb{Q}$. Moreover, we see that that $pqr \equiv 1 \pmod {4}$ and consequently, $2$ is not totally ramified in $F/\mathbb{Q}$. Thus the ramification index of the rational prime $2$ in $F/\mathbb{Q}$ is also $2$. Therefore, by Proposition \ref{zantema's-main-theorem}, it follows that $H^{1}(G,\mathcal{O}_{F}^{*})$ injects into $\mathbb{Z}/2\mathbb{Z} \oplus \mathbb{Z}/2\mathbb{Z} \oplus \mathbb{Z}/2\mathbb{Z} \oplus \mathbb{Z}/2\mathbb{Z}$. Since the rational prime $2$ is not totally ramified in $F/\mathbb{Q}$, by Lemma \ref{lem-1}, we have $H = H^{1}(G,\mathcal{O}_{F}^{*})$. We prove that $|H| = 8$. This, together with Proposition \ref{zantema's-main-theorem}, yields that $Po(F) \simeq \mathbb{Z}/2\mathbb{Z}$.

\smallskip

In the notations of Lemma \ref{lem-2}, we have $\Delta_{1} = p$, $\Delta_{2} = qr$ and $\Delta_{3} = pqr$. Also, since $p \equiv 3 \pmod {4}$, a fundamental unit $u_{1} = z_{1} + t_{1}\sqrt{p}$ of $F_{1}$ necessarily has norm $1$. In the notations of Lemma \ref{lem-0}, we have $$a_{1} = N(u_{1} + 1) = (z_{1} + 1)^{2} - pt_{1}^{2} = 2(z_{1} + 1).$$ Consequently, by Lemma \ref{lem-0}, we conclude that $[a_{1}] = [2(z_{1} + 1)] = [2] \mbox{ or } [2p]$. Thus we obtain that $$[a_{1}], [\Delta_{1}], [\Delta_{2}], [\Delta_{3}] \in \langle [2],[p],[qr] \rangle \subseteq \mathbb{Q}^{*}/(\mathbb{Q}^{*})^{2}.$$ Since $\left(\frac{q}{r}\right) = -1$, we conclude by Dirichlet's theorem (cf. \cite{trotter}) that the fundamental unit of $F_{2}$ has norm $-1$ and therefore, $[a_{2}] = [1]$. 

\smallskip

Now, let $u_{3} = \frac{1}{2}(z + t\sqrt{pqr})$ be a fundamental unit of $F_{3}$. Again, since $p \equiv 3 \pmod {4}$, we cannot have $N(u_{3}) = -1$. Therefore, $N(u_{3}) = \frac{z^{2} - t^{2}pqr}{4} = 1$. That is, $(z - 2)(z + 2) = t^{2}pqr$. We also notice that $\gcd(z - 2, z + 2) = 1 \mbox{ or } 2 \mbox{ or } 4$. 

\smallskip

\noindent
{\bf Case 1.} $\gcd(z - 2, z + 2) = 1$. Then from the equation $(z - 2)(z + 2) = t^{2}pqr$, we have $z - 2 = n^{2}\eta$ and $z + 2 = m^{2}\epsilon$, where $m,n,\eta$ and $\epsilon$ are integers with $\gcd(m,n) = 1$, $mn = t$ and $\eta\epsilon = pqr$.

\smallskip

Now, using the equation $z^{2} - t^{2}pqr = 4$, we have $a_{3} = N(u_{3} + 1) = z + 2$. Hence $[a_{3}] = [z + 2] = [m^{2}\epsilon] = [\epsilon] \in \mathbb{Q}^{*}/(\mathbb{Q}^{*})^{2}$. Thus $[z + 2] \in \langle [2],[p],[qr] \rangle$ if and only if $\epsilon = 1,p,qr,pqr$. We prove that $\epsilon = q,r,pq,pr$ cannot happen. This proves that $a_{3} \in \langle [2],[p],[qr] \rangle$.

\smallskip

If $\epsilon = q$, then $\eta = pr$ and therefore, $$4 = (z + 2) - (z - 2) = m^{2}\epsilon - n^{2}\eta = m^{2}q - n^{2}pr.$$ Consequently, $1 = \left(\frac{4}{r}\right) = \left(\frac{q}{r}\right) = -1$ yields a contradiction. Again, if $\epsilon = r$, then $\eta = pq$ and $4 = m^{2}r - n^{2}pq$. Thus $1 = \left(\frac{4}{q}\right) = \left(\frac{r}{q}\right) = -1$ yields a contradiction. Similarly, for the cases when $\epsilon = pq$ or $\epsilon = pr$, reading the resulting equations modulo $r$ leads to a contradiction.


\medskip

\noindent
{\bf Case 2.} $\gcd(z - 2, z + 2) = 2$. Then $t$ is an even integer and hence $t = 2k$ for some integer $k$. Therefore, the equation $\frac{z - 2}{2}\cdot \frac{z + 2}{2} = k^{2}pqr$, together with the fact that $\gcd(\frac{z - 2}{2}, \frac{z + 2}{2}) = 1$, implies that $\frac{z - 2}{2} = n^{2}\eta$ and $\frac{z + 2}{2} = m^{2}\epsilon$ for integers $m,n,\eta$ and $\epsilon$ such that $\gcd(m,n) = 1$, $mn = k$ and $\eta\epsilon = pqr$. Consequently, we have $[a_{3}] = [z + 2] = [2m^{2}\epsilon] = [2\epsilon] \in \langle [2],[p],[qr] \rangle$ if and only if $\epsilon = 1,p,qr,pqr$. We now prove that $\epsilon = q,r,pq,pr$ cannot happen.

\smallskip

If $\epsilon = q$, then $\eta = pr$ and therefore, $$2 = \frac{z + 2}{2} - \frac{z - 2}{2} = m^{2}\epsilon - n^{2}\eta = m^{2}q - n^{2}pr.$$ Consequently, $1 = \left(\frac{2}{r}\right) = \left(\frac{m^{2}q}{r}\right) = \left(\frac{q}{r}\right) = -1$ yields a contradiction. Again, if $\epsilon = r \mbox{ or } pq \mbox{ or } pr$, then reading the resulting equations modulo $r$ and using the hypotheses on the Legendre symbols yield a contradiction.


\medskip

\noindent
{\bf Case 3.} $\gcd(z - 2, z+ 2) = 4$. Then $t = 4k$ for some integer $k$ and therefore, $\frac{z - 2}{4}\cdot \frac{z + 2}{4} = k^{2}pqr$. Since $\gcd(\frac{z - 2}{4}, \frac{z + 2}{4}) = 1$, we have $\frac{z - 2}{4} = n^{2}\eta$ and $\frac{z + 2}{4} = m^{2}\epsilon$ for integers $m,n,\eta, \epsilon$ with $\gcd(m,n) = 1$, $mn = k$ and $\epsilon\eta = pqr$. Consequently, we have $[a_{3}] = [z + 2] = [4m^{2}\epsilon] = [\epsilon] \in \langle [2],[p],[qr] \rangle$ if and only if $\epsilon = 1,p,qr,pqr$. We now prove that $\epsilon = q,r,pq,pr$ cannot happen.

\smallskip

If $\epsilon = q$, then $\eta = pr$ and therefore, $$1 = \frac{z + 2}{4} - \frac{z - 2}{4} = m^{2}\epsilon - n^{2}\eta = m^{2}q - n^{2}pr.$$ Consequently, $1 = \left(\frac{1}{r}\right) = \left(\frac{m^{2}q}{r}\right) = \left(\frac{q}{r}\right) = -1$ yields a contradiction. Again, if $\epsilon = r$, then $\eta = pq$ and $1 = m^{2}r - n^{2}pq$. Therefore, $1 = \left(\frac{1}{q}\right) = \left(\frac{m^{2}r}{q}\right) = \left(\frac{r}{q}\right) = -1$ yields a contradiction. Similarly, if $\epsilon = pq \mbox{ or } pr$, then considering the equation $1 = m^{2}\epsilon - n^{2}\eta$ modulo $r$ leads to a contradiction.


\smallskip

Thus, combining all the above cases, we conclude that $H = \langle [2],[p],[qr] \rangle \subseteq \mathbb{Q}^{*}/(\mathbb{Q}^{*})^{2}$. Hence $|H^{1}(G,\mathcal{O}_{F}^{*})| = |H| = 8$ and therefore by Lemma \ref{zantema's-main-theorem}, we have $Po(F) \simeq \mathbb{Z}/2\mathbb{Z}$. This completes the proof of Theorem \ref{NEW-TH}. $\hfill\Box$

\medskip

\section{Proof of Theorem \ref{main-2}}

Let $K_{1} = \mathbb{Q}(\sqrt{2})$, $K_{2} = \mathbb{Q}(\sqrt{pq})$ and $K_{3} = \mathbb{Q}(\sqrt{2pq})$. Then the rational primes that are ramified in $K^{\prime \prime} = \mathbb{Q}(\sqrt{2},\sqrt{pq})$ are precisely $2,p$ and $q$. Moreover, the ramification indices are each equal to $2$. Therefore, by Proposition \ref{zantema's-main-theorem}, it follows that $H^{1}(G,\mathcal{O}_{K^{\prime \prime}}^{*})$ injects into $\mathbb{Z}/2\mathbb{Z} \oplus \mathbb{Z}/2\mathbb{Z} \oplus \mathbb{Z}/2\mathbb{Z}$. Since the rational prime $2$ is not totally ramified in $K^{\prime \prime}/\mathbb{Q}$, by Lemma \ref{lem-1}, we have $H = H^{1}(G,\mathcal{O}_{K^{\prime \prime}}^{*})$. We prove that $|H| = 4$. Then by Proposition \ref{zantema's-main-theorem}, it follows that $Po(K^{\prime \prime}) \simeq \mathbb{Z}/2\mathbb{Z}$.

\smallskip

Now, in the notations of Lemma \ref{lem-2}, we have $\Delta_{1} = 2, \Delta_{2} = pq, \Delta_{3} = 2pq$. Also, since $1 + \sqrt{2}$ is a fundamental unit of $\mathbb{Q}(\sqrt{2})$, we see that $N(1 + \sqrt{2}) = -1$ and therefore, $a_{1} = 1$. Again, since $\left(\frac{p}{q}\right) = -1$, by Dirichlet's theorem (also cf. \cite{trotter}), we conclude that the fundamental unit of $K_{2} = \mathbb{Q}(\sqrt{pq})$ has norm $-1$. Therefore, $a_{2} = 1$.

\smallskip

Now, let $u = z + t\sqrt{2pq}$ be a fundamental unit of $K_{3} = \mathbb{Q}(\sqrt{2pq})$. If $N(u) = -1$, then $a_{3} = 1$. In that case, we have $$\langle [\Delta_{1}], [\Delta_{2}], [\Delta_{3}], [a_{1}], [a_{2}], [a_{3}]\rangle = \langle [2], [pq] \rangle \subseteq \mathbb{Q}^{*}/(\mathbb{Q}^{*})^{2}$$ and therefore, $|H| = 4$. Otherwise, $N(u) = z^{2} - 2t^{2}pq = 1$. That is $(z-1)(z+1) = 2t^{2}pq$. Since $(z + 1) - (z - 1) = 2$, we conclude that $\gcd(z - 1, z + 1) = 1 \mbox{ or } 2$. Also, $2$ divides $2t^{2}pq$ and consequently $z$ is odd. Thus we have $\gcd(z - 1, z + 1) = 2$. Therefore, we have $t = 2k$ for some integer $k$ and hence 
\begin{equation}\label{NEW-Equn}
\frac{z - 1}{2}\cdot \frac{z + 1}{2} = 2k^{2}pq.
\end{equation}

Now, the integers $\frac{z - 1}{2}$ and $\frac{z + 1}{2}$ are relatively prime and therefore from equation \eqref{NEW-Equn}, we conclude that $\frac{z - 1}{2} = n^{2}\eta$ and $\frac{z + 1}{2} = m^{2}\epsilon$, where $m, n$ are coprime integers such that $mn=t$ and $\epsilon, \eta$ are coprime integers such that $\epsilon\eta = 2pq$.

\smallskip

Therefore, we have $[2(z+1)] = [4m^{2}\epsilon] = [\epsilon] \in \mathbb{Q}^{*}/(\mathbb{Q}^{*})^{2}$. Since $\epsilon$ is a divisor of $2pq$, we have $\epsilon \in \{1,2,p,q,2p,2q,pq,2pq\}$. We immediately see that if $\epsilon = 1 \mbox{ or } 2 \mbox{ or } pq \mbox{ or } 2pq$, then $[\epsilon] \in \langle [2],[pq] \rangle$. We prove that the remaining choices for $\epsilon$ cannot occur.

\smallskip

We first note that eliminating $z$ from the equations $\frac{z+1}{2} = m^{2}\epsilon$ and $\frac{z-1}{2} = n^{2}\eta$, we obtain $1 = m^{2}\epsilon - n^{2}\eta$. Now, for $\epsilon = p$, we have $\eta = 2q$ and consequently, $1 = m^{2}p - 2n^{2}q$. Hence $1 = \left(\frac{1}{q}\right) = \left(\frac{p}{q}\right) = -1$, which is a contradiction. Similarly, for $\epsilon = q$, we have $\eta = 2p$ and therefore, $1 = m^{2}q - 2n^{2}p$. Hence $1 = \left(\frac{1}{p}\right) = \left(\frac{q}{p}\right) = -1$, which is a contradiction. Again for $\epsilon = 2p$, we have $\eta = q$ and therefore, $1 = 2m^{2}p - n^{2}q$. Hence $1 = \left(\frac{1}{p}\right) = \left(\frac{q}{p}\right) = -1$, which yields a contradiction. Lastly, for $\epsilon = 2q$, we have $\eta = p$ and therefore, $1 = 2m^{2}q - n^{2}p$. Hence $1 = \left(\frac{1}{q}\right) = \left(\frac{p}{q}\right) = -1$, which is a contradiction.

\smallskip

Thus we obtain that $$\langle [\Delta_{1}], [\Delta_{2}], [\Delta_{3}], [a_{1}], [a_{2}], [a_{3}]\rangle = \langle [2], [pq] \rangle \subseteq \mathbb{Q}^{*}/(\mathbb{Q}^{*})^{2}.$$ 

Consequently, we have $|H^{1}(G,\mathcal{O}^{*}_{K^{\prime \prime}})| = |H| = 4$ and hence $Po(K^{\prime \prime}) \simeq \mathbb{Z}/2\mathbb{Z}$. This completes the proof of Theorem \ref{main-2}. $\hfill\Box$


\section{Concluding remarks}
Theorem \ref{main-1} and Theorem \ref{main-2} provide new families of bi-quadratic number fields with P\'{o}lya groups isomorphic to $\mathbb{Z}/2\mathbb{Z}$. From our previous discussion on Euclidean ideal class in Section \ref{subsection}, we can link the two interesting features of a number field, namely, having an Euclidean ideal class and being a P\'{o}lya field. In \cite{self-jnt}, two lists of bi-quadratic fields with an EIC were given where each of the number fields $\mathbb{Q}(\sqrt{p},\sqrt{qr})$ and $\mathbb{Q}(\sqrt{2},\sqrt{qr})$ has class number $2$.

We have observed that plenty of number fields from those list satisfy the hypotheses of Theorem \ref{main-1} and Theorem \ref{main-2}. For Theorem \ref{main-1}, we see that $\mathbb{Q}(\sqrt{3},\sqrt{17\cdot 29})$, $\mathbb{Q}(\sqrt{3},\sqrt{29\cdot 41})$, $\mathbb{Q}(\sqrt{3},\sqrt{29\cdot 113})$ etc. are admissible choices. For Theorem \ref{main-2}, we provide a list as follows.


\begin{table}[h!]
	\centering
	\begin{tabular}{||c c c c ||} 
		\hline
		$(2,p,q)$ & $(2,p,q)$ & $(2,p,q)$ & $(2,p,q)$\\ [0.5ex]
		\hline\hline
(2, 5, 17) & (2, 5, 37) & (2, 5, 97) & (2, 5, 173)\\
(2, 5, 193) & (2, 13, 37) & (2, 13, 73) & (2, 13, 89)\\ 
(2, 13, 97) & (2, 13, 109) & (2, 13, 193) & (2, 13, 197)\\ 
(2, 17, 5) & (2, 17, 29) & (2, 17, 37) & (2, 17, 61)\\
(2, 17, 197) & (2, 29, 17) & (2, 29, 61) & (2, 29, 89)\\  
 [1ex]
\hline
\end{tabular}
\end{table}

\smallskip

In fact, we proved in this paper that even though some number fields have a non-principal EIC, but they are not P\'{o}lya fields. In the light of this, we end this section with the two following questions.

\begin{question}
Let $K$ be a number field with an Euclidean ideal class. When is $K$ a P\'{o}lya field?
\end{question}

\begin{question}
Let $K$ be a number field with non-trivial P\'{o}lya group $Po(K)$. Assume that $Po(K)$ is cyclic. When is $Cl_{K}$ cyclic?
\end{question}

\medskip

{\bf Acknowledgements.} The first author gratefully acknowledges the financial support provided by Indian Institute of Technology, Guwahati to carry out this research.

\end{document}